\documentclass[12pt,a4paper,oneside]{amsart}

\usepackage{amssymb} 
\usepackage{amsmath}  
\usepackage{amsthm}
\usepackage{amsfonts}

\usepackage[english]{babel}
\usepackage[utf8]{inputenc}
\usepackage{enumerate}

\usepackage{mathrsfs}  
\usepackage{graphicx}

\usepackage{float} 

\usepackage[inline]{enumitem}

\usepackage[colorinlistoftodos,prependcaption,textsize=tiny]{todonotes}

\vbadness=\maxdimen

\newtheorem{thm}{Theorem}[section]
\newtheorem*{thm*}{Theorem}
\newtheorem{prop}[thm]{Proposition}
\newtheorem{lem}[thm]{Lemma}
\newtheorem{rem}[thm]{Remark}
\newtheorem{cor}[thm]{Corollary}
\newtheorem{obs}[thm]{Observation}
\newtheorem{df}[thm]{Definition}

\newtheorem{fact}[thm]{Fact}

\newcommand{\mm}{\mathscr {M} _{m}} 
\newcommand{\Pm}{\mathscr {P} _{m}} 
\newcommand{\Pmm}{\mathscr {P} _{m,m}} 

\newcommand{\bdf}{\begin{df}} 
\newcommand{\edf}{\end{df}}
\newcommand{\bthm}{\begin{thm}} 
\newcommand{\ethm}{\end{thm}}
\newcommand{\bcor}{\begin{cor}} 
\newcommand{\ecor}{\end{cor}}
\newcommand{\blem}{\begin{lem}} 
\newcommand{\elem}{\end{lem}}

\begin{document}

\title{Geodesics in a Graph of Perfect Matchings}
\author{Roy H. Jennings}
\address{Department of Mathematics, Bar-Ilan University, Ramat-Gan 52900, Israel.}
\email{RoyHJennings@gmail.com}
\date{}

\begin{abstract}
\footnotesize
Let $\mathscr{P}_{m}$ be the graph on the set of perfect matchings in the complete graph $K_{2m}$, where two perfect matchings are connected by an edge if their symmetric difference is a cycle of length four. 
This paper studies geodesics in $\mathscr{P}_{m}$. 
The diameter of $\mathscr{P}_{m}$, as well as the eccentricity of each vertex, are shown to be $m-1$. 
Two proof are given to show that the number of geodesics between any two antipodes is $m^{m-2}$. The first is a direct proof via a recursive formula, and the second is via reduction to the number of minimal factorizations of a given $m$-cycle in the symmetric group $S_m$.
An explicit formula for the number of geodesics between any two matchings in $\mathscr{P}_{m}$ is also given.

Let $\mathscr{M}_m$ be the graph on the set of non-crossing perfect matchings of $2m$ labeled points on a circle with the same adjacency condition as in $\mathscr{P}_m$. $\mathscr{M}_m$ is an induced subgraph of $\mathscr{P}_m$, and it is shown that $\mathscr{M}_m$ has exactly one pair of antipodes having the maximal number ($m^{m-2}$) of geodesics between them.
\end{abstract}

\maketitle

\setcounter{page}{1}
\section{Introduction}
Consider a set of $2m$ labeled points on a circle. Join its points in disjoint pairs by $m$ straight line segments, such that no two lines intersect.
Abstractly, this is a perfect matching in the complete graph $K_{2m}$. 
Such a matching is called a {\bf non-crossing perfect matching}. 
Hernando, Hurtado and Noy defined in \cite{noy} the {\bf graph of non-crossing perfect matchings} $\mm$, in which two matchings are connected by an edge if their symmetric difference is a cycle of length four. They showed that the diameter and the eccentricity of every vertex in this graph are $m-1$.

This paper studies the graph $\Pm$ on the set of all perfect matchings in the complete graph $K_{2m}$ where two matchings are connected by an edge if their symmetric difference is a cycle of length four.

The distance between any two matchings in $\Pm$ is shown to depend only on the number of components in their union.

\begin{thm*}
[Theorem \ref{thm_distance}]
For any $M_1, M_2 \in \Pm$, $d(M_1,M_2) = m-l$, where $l$ is the number of components in $M_1 \cup M_2$.
\end{thm*}
Although $\Pm$ is larger than the graph of non-crossing perfect matchings $\mm$, which is an induced subgraph of $\Pm$, it is shown that the two graphs still share the same diameter, and that it is equal to the eccentricity of all of their vertices.
\begin{thm*}
[Corollary \ref{cor_diam_ecc}] The diameter of $\Pm$, as well as the eccentricity of every vertex in it, is $m-1$.
\end{thm*}
Enumeration of the number of geodesics (shortest paths) between antipodes in $\Pm$ reveals the following surprising result.
\begin{thm*}
[Corollary \ref{cor_antipodes_numer_of_paths}] The number of geodesics between any two antipodes in $\Pm$ is $m^{m-2}$.
\end{thm*}
The expression $m^{m-2}$ appears in Cayley's well known formula for the number of labeled trees on $m$ vertices. It was also proved by D\'{e}nes \cite{denes} (see also \cite{hurwitz, strehl}) to be equal to the number of factorizations of a given $m$-cycle as a product of $m-1$ transpositions in $S_m$. 
In Section \ref{section_cycleFactorization}, an alternative proof, suggested by Y. Roichman, for Corollary \ref{cor_antipodes_numer_of_paths} is provided via reduction to the number of factorizations of a given $m$-cycle as a product of $m-1$ transpositions. 

An explicit formula for the number of geodesics between any two matchings in $\Pm$ is given in Corollary \ref{cor_number_of_geodesics}.

Corollary \ref{cor_antipodes_numer_of_paths} does not hold for the subgraph $\mm$. In fact
\begin{thm*}
[Theorem \ref{thm_non_cross_antipods}]
The graph $\mm$ of non-crossing perfect matchings has a unique pair of matchings with $m^{m-2}$ geodesics between them.
All other pairs have a smaller number of geodesics.
\end{thm*}

This paper is based on the author's M.Sc.\ thesis. For generalizations and extensions of results in this paper see \cite{roichmann, avni, cohen}.

\section{The Graph of Perfect Matchings}
Consider the set of all perfect matchings in the complete graph $K_{2m}$. Denote by $\Pm$ the \textbf{graph of perfect matchings} on this set, in which two distinct matchings $M_1$ and $M_2$ are connected by an edge, if their symmetric difference is a cycle of length four. Here, the symmetric difference of $M_1$ and $M_2$, is the graph consisting of the edges that belong to exactly one of these matchings. In this work all matchings are perfect, and thus the term {\it perfect} is freely omitted.
Denote adjacent matchings in $\Pm$ by $M_1\sim M_2$, and write $M_1 \simeq M_2$ if either $M_1=M_2$ or $M_1 \sim M_2$.

\begin{figure}[hbt]
\[
\begin{aligned}
\begin{tikzpicture}[scale=0.5]
\fill (2.4,3.4) circle (0.1); 
\fill (3.4,2.4) circle (0.1); 
\fill (3.4,1)   circle (0.1); 
\fill (2.4,0)   circle (0.1); 
\fill (1,0)     circle (0.1); 
\fill (0,1)     circle (0.1); 
\fill (0,2.4)   circle (0.1); 
\fill (1,3.4)   circle (0.1); 

\draw (0,1)-- (0,2.4);
\draw (1,3.4)-- (1,0);
\draw[red] (2.4,3.4)-- (3.4,2.4);
\draw[red] (3.4,1)-- (2.4,0);
\draw[blue,dotted] (2.4,3.4)-- (2.4,0);
\draw[blue,dotted] (3.4,1)-- (3.4,2.4);
\end{tikzpicture}
\end{aligned}
\ \ \ \longleftrightarrow \ \ \
\begin{aligned}
\begin{tikzpicture}[scale=0.5]
\fill (2.4,3.4) circle (0.1); 
\fill (3.4,2.4) circle (0.1); 
\fill (3.4,1)   circle (0.1); 
\fill (2.4,0)   circle (0.1); 
\fill (1,0)     circle (0.1); 
\fill (0,1)     circle (0.1); 
\fill (0,2.4)   circle (0.1); 
\fill (1,3.4)   circle (0.1); 

\draw (0,1)-- (0,2.4);
\draw (1,3.4)-- (1,0);
\draw[red,dotted] (2.4,3.4)-- (3.4,2.4);
\draw[red,dotted] (3.4,1)-- (2.4,0);
\draw[blue] (2.4,3.4)-- (2.4,0);
\draw[blue] (3.4,1)-- (3.4,2.4);
\end{tikzpicture}
\end{aligned}
\]
\caption{Adjacent matchings in $\mathscr{P}_4$}
\label{fig:ex_match_flip}
\end{figure}

\begin{fact}
The cardinality of $\Pm$ is the double factorial:
$$|\Pm|=(2m-1)!!=(2m-1)\cdot (2m-3) \cdot ... \cdot 1.$$
\end{fact}

For any edge $e=(v_1,v_2)$ in the complete graph $K_{2m}$ and a matching $M$ in $\Pm$, the \textbf{insertion of the edge $e$ into the matching $M$}, denoted by $M*e$, is defined as follows:\\
If $e$ is already in $M$, then $M*e = M$. Otherwise, if $(v_1,v_3)$ and $(v_2,v_4)$ are the edges in $M$ incident with $v_1$ and $v_2$, then $M*(v_1,v_2)$ is the matching obtained from $M$ by deleting $(v_1,v_3)$ and $(v_2,v_4)$ and adding $(v_1,v_2)$ and $(v_3,v_4)$.

Extend the definition recursively:
$$ M*(e_1,\dots,e_n) = (M*(e_1,\dots,e_{n-1}))*e_n.$$

\begin{obs}
Two distinct matchings $M_1$ and $M_2$ in $\Pm$ are adjacent if and only if $M_2 = M_1*e$, for some edge $e$ not in $M_1$. 
\end{obs}

\begin{rem} \label{rem_regular}
Let $M$ be a matching in $\Pm$. By definition, every neighbor of $M$ misses exactly two edges of $M$. Conversely, for each pair of distinct edges in $M$, there are exactly two neighbors of $M$ not containing this pair. Therefore, $\Pm$ is a regular graph of degree $2 \binom{m}{2}$.
\end{rem}

\begin{prop}
Let $M$ be a matching in $\Pm$, and let $(e_1, \dots, e_m)$ be the edges of $M$ arranged in any order. For any matching $M'$ in $\Pm$, $M = M' * (e_1, \dots, e_m)$.
\end{prop}
Indeed, the edges inserted into $M'$ are vertex disjoint, and therefore they all belong
to the final matching.

\begin{cor}
\label{cor_cm_connected}
The graph of perfect matchings $\Pm$ is connected.
\end{cor}

\section{Geodesics and Diameter}
Recall that a {\bf geodesic} between two vertices in a connected graph, is a shortest path between them.
The {\bf distance} between two vertices $u$ and $v$, denoted by $d(u,v)$, is the length of a geodesic between them.
The {\bf eccentricity} of a vertex is the maximal distance between this vertex and any other vertex.
The maximal length of a geodesic in a connected graph is called the graph's {\bf diameter}.
If the distance between two vertices is the graph's diameter, then the vertices are called {\bf antipodes}.

Hernando, Hurtado and Noy \cite{noy} found a formula for the distance between two matchings in the graph $\mm$ of non-crossing perfect matchings with $m$ edges, implying that the diameter of this graph and the eccentricity of every vertex in it are $m-1$. In this section, it is shown that although $\Pm$ is larger than the graph of non-crossing perfect matchings, where the latter is an induced subgraph, the diameter and the eccentricity of every vertex in it remain $m-1$.

\begin{fact}
\label{fact_diag_union}
The union of two matchings $M_1$ and $M_2$ in $\Pm$, is a vertex-disjoint union of alternating cycles $C_i$ of even length (a common edge is considered as a cycle of length two):
$$M_1 \cup M_2 = C_1\;\dot\cup\; C_2 \;\dot\cup\;\dots \;\dot\cup\; C_l\ .$$
\end{fact}

\begin{obs}
\label{obs_inc_comp}
Denote the number of connected components in a graph $G$ by $c(G)$.
Given two matchings $M_1, M_2 \in \Pm$ and an edge $e \notin M_1$, exactly one of the following cases holds:
\begin{enumerate}
\item $c(M_1*e\cup M_2)=c(M_1\cup M_2)-1$. This is the case if and only if the vertices of $e$ belong to two different cycles in $M_1 \cup M_2$.
\item\label{obs_inc_comp_split} $c(M_1*e\cup M_2)=c(M_1\cup M_2)+1$. This is the case if and only if the vertices of $e$ belong to the same cycle in $M_1 \cup M_2$, and the insertion of $e$ into $M_1$ splits this cycle into two cycles in $M_1*e \cup M_2$.
\item $c(M_1*e\cup M_2)=c(M_1\cup M_2)$. This is the case if and only if the vertices of $e$ belong to the same cycle in $M_1 \cup M_2$, and the insertion of $e$ into $M_1$ does not split this cycle in $M_1*e \cup M_2$.
\end{enumerate}
\end{obs}

\begin{rem}\label{rem_inc_comp}
In the settings of Observation \ref{obs_inc_comp}, every pair of distinct edges in $M_1$ that belong to different components in $M_1\cup M_2$, corresponds to two unique neighbors of $M_1$ that belong to the first case of the observation. 
And every pair of distinct edges in $M_1$ that belongs to the same component in $M_1\cup M_2$, corresponds to one unique neighbor of $M_1$ that belongs to the second case of the observation, and one unique neighbor that belongs to the third case.
\end{rem}

\begin{thm}
\label{thm_distance}
For any $M_1, M_2 \in \Pm$, $d(M_1,M_2) = m-l$, where $l$ is the number of components in $M_1 \cup M_2$.
\end{thm}
\begin{proof}
Note that $m-l=0$ if and only if $M_1=M_2$.
By Observation \ref{obs_inc_comp}, $\left|c(M'\cup M_2)-c(M\cup M_2)\right|\leq 1$ for any two neighbors $M\sim M'$ in $\Pm$.
Thus $d(M_1,M_2)\geq m-l$.
For every matching $M\neq M_2$, there is some alternating cycle $C\subseteq M\cup M_2$ such that $|C|\geq 4$. Therefore, by Remark \ref{rem_inc_comp}, for every pair of distinct edges in $M\cap C$, $M$ has a unique neighbor $M'$ such that $c(M'\cup M_2)-c(M\cup M_2) = 1$.
Thus, $d(M_1,M_2)\leq m-l$.
\end{proof}

By Theorem \ref{thm_distance}, antipodes in $\Pm$ are pairs of matchings whose union consists of one cycle.
\begin{cor}
\label{cor_diam_ecc}
The diameter of $\Pm$, as well as the eccentricity of every vertex in it, are $m-1$.
The number of antipodes of every matching is $(2m-2)!!$.
\end{cor}

\begin{rem}\label{rem_noy}
An analogue of Theorem \ref{thm_distance} for the graph $\mm$ of non-crossing perfect matchings appears in \cite{noy}, in the equivalent form $$d(M_1,M_2) = \frac{1}{2} \sum_{i=1}^{l} \left(\operatorname{length}(C_i) - 2\right).$$
\end{rem}

\section{Counting Geodesics}
In this section the number of geodesics between any two matchings in $\Pm$ is given. In particular, the number of geodesics between antipodes is shown to be $m^{m-2}$. The section concludes by showing that in the subgraph $\mm$ of non-crossing perfect matchings, there is exactly one pair of antipodes having a maximal number ($m^{m-2}$) of geodesics between them.

\bdf
Denote by $P_{2k} \ (k\geq 2)$ the number of geodesics between two matchings, whose symmetric difference is one cycle of length $2k$. Define $P_2 = 1$.
\edf



\bthm For every positive integer $k$, we have
\label{thm_p2k}
$$
P_{2k} = \frac{k}{2}\sum_{i=1}^{k-1} {\binom{k-2}{i-1}} P_{2i}P_{2k-2i} \ \ (k\geq 2),
$$
with $P_2=1$.
\ethm

\begin{proof}
Let $M_1, M_2 \in \Pm$ ($m\geq k$) be two matchings with a symmetric difference of one cycle of length $2k$. Denote this cycle by $C$. $M_1 \cup M_2$ has $m-k+1$ components ($m-k$ components of size 2 and one component of size $2k$). 
Therefore, by Theorem \ref{thm_distance}, $d(M_1, M_2)=k-1$. 

By Remark \ref{rem_inc_comp} and Theorem \ref{thm_distance}, a neighbor $M'$ of $M_1$ is in some geodesic between $M_1$ and $M_2$ if and only if $C$ splits in $M'\cup M_2$ into two cycles. 
The lengths of these two cycles are $2i$ and $2k-2i$ for some $1\leq i\leq\frac{k}{2}$. 
Therefore, the number of geodesics from $M_1$ to $M_2$, beginning with $M'$, is ${\binom{k-2}{i-1}}P_{2i}P_{2k-2i}$, where the binomial coefficient ${\binom{k-2}{i-1}}$ counts the ways to interlace the remaining $k-2$ insertions between the two cycles.

For every $1\leq i< \frac{k}{2}$, there are exactly $k$ neighbors of $M_1$ in which $C$ splits into two cycles of length $2i$ and $2k-2i$ as above. If $k$ is even, then for $i=\frac{k}{2}$ there are exactly $\frac{k}{2}$ neighbors of $M_1$ in which $C$ splits into two cycles of length $k$ as above. Therefore, the formula $k\sum_{i=1}^{k-1} {\binom{k-2}{i-1}} P_{2i}P_{2k-2i}$ counts every geodesic twice, leading to the claimed result.
\end{proof}

The following is a reformulation of Theorem \ref{thm_p2k}, and its proof was suggested by R. Adin.
\begin{lem} For every positive integer $k$, we have
\label{lem_t2k_p2k}
$$
P_{2k} = \sum_{i=1}^{k-1} i{\binom{k-2}{i-1}} P_{2i}P_{2k-2i},
$$
with $P_2=1$.
\end{lem}

\begin{proof} Denote $a_i = i{\binom{k-2}{i-1}} P_{2i}P_{2k-2i}$. The symmetry of the binomial coefficients ${\binom{k-2}{i-1}} = {\binom{k-2}{k-i-1}}$ implies 

$ \begin{array}{lll}
2\sum_{i=1}^{k-1} a_i & = & \sum_{i=1}^{k-1}a_i + \sum_{i=1}^{k-1} a_{k-i} \\
  & = & \sum_{i=1}^{k-1}\left( i{\binom{k-2}{i-1}} P_{2i}P_{2k-2i} + (k-i){\binom{k-2}{k-i-1}}P_{2(k-i)}P_{2k-2(k-i)}\right) \\
  & = & \sum_{i=1}^{k-1}\left( i{\binom{k-2}{i-1}} P_{2i}P_{2k-2i} + (k-i){\binom{k-2}{i-1}}P_{2k-2i}P_{2i}\right) \\
  & = & k \sum_{i=1}^{k-1} {\binom{k-2}{i-1}} P_{2i}P_{2k-2i} \\
  & = & 2P_{2k}
\end{array} $

The last equation is essentially Theorem \ref{thm_p2k}.
\end{proof}

\begin{cor}
\label{cor_numer_of_paths}
For every positive integer $k$, we have $P_{2k} = k^{k-2}.$
\end{cor}

\begin{proof}
By the well known Cayley Formula \cite{cayley}, the number of labeled trees on $k$ vertices is $T_k = k^{k-2}$.
By \cite[ex. 6, p. 34 and pp. 249-250]{lovaasz}, we have $T_k = \sum_{i=1}^{k-1} i{\binom{k-2}{i-1}} T_{i}T_{k-i}$, with $T_1 = 1$. 
Comparison with Lemma \ref{lem_t2k_p2k} concludes the proof.
\end{proof}

By Theorem \ref{thm_distance}, antipodes in $\Pm$ are pairs of matchings with a symmetric difference of one cycle of length $2m$.

\begin{cor}
\label{cor_antipodes_numer_of_paths}
The number of geodesics between antipodes in $\Pm$ is $m^{m-2}$.
\end{cor}

Corollary \ref{cor_antipodes_numer_of_paths} can be generalized to count the number of geodesics between any two matchings in $\Pm$.

\bcor\label{cor_number_of_geodesics}
Let $M_1, M_2 \in \Pm$ with $M_1 \cup M_2 = C_1 \;\dot\cup\; C_2 \;\dot\cup\; \dots \;\dot\cup\; C_l$, as in Fact \ref{fact_diag_union}. 
The number of geodesics between $M_1$ and $M_2$ is
$$\binom{m-l}{n(C_1),\dots,n(C_l)} \prod_{i=1}^{l}\left(n(C_i)+1\right)^{n(C_i)-1},$$
where $n(C_i)=\frac{length(C_i)}{2} - 1$.
\ecor

\begin{proof}
$n(C_i)$ is the number of insertions within the cycle $C_i$. The binomial coefficient ${\binom{m-l}{n(C_1),\dots,n(C_l)}}$ counts the ways to interlace insertions between the cycles.
\end{proof}

The analogue of Corollary \ref{cor_antipodes_numer_of_paths} for the induced subgraph in $\Pm$ of non-crossing perfect matchings $\mm$, is as follows.

\begin{thm} \label{thm_non_cross_antipods}
The graph $\mm$ of non-crossing perfect matchings has a unique pair of matchings with $m^{m-2}$ geodesics between them.
All other pairs have a smaller number of geodesics.
\end{thm}

\begin{proof}
$\mm$ consists of matchings on a set of $2m$ points on a circle. Denote the convex hull of these points by $H$.
We show that the two matchings having all their edges on the boundary of $H$, are the only pair of antipodes in $\mm$ with $m^{m-2}$ geodesics between them.
All other pairs have a smaller number of geodesics.

Let $M_1$ and $M_2$ be the two matchings in $\mm$ having all their edges on the boundary of $H$. By Remark \ref{rem_noy}, $M_1$ and $M_2$ are antipodes in $\mm$. By Theorem \ref{thm_distance}, they are also antipodes in $\Pm$. 
Let $P=(M_1=M'_1,M'_2,\dots,M'_d=M2)$ be a geodesic between $M_1$ and $M_2$.
Since $M'_1\cup M'_d$ has one component, by Theorem \ref{thm_distance} $M'_2 \cup M'_d$ is a union of two vertex-disjoint alternating cycles. 
It is clear that these cycles are the boundaries of two disjoint polytopes (including the case of a shared edge, which is a convex polytope with two vertices).
Therefore, $M'_2$ is also in $\mm$.
Similarly, for every $M'_i\in P$, $2\leq i\leq d$, some convex polytope in $M'_{i-1}\cup M'_d$, is split into two disjoint convex polytopes in $M'_{i}\cup M'_d$.
Thus, $P$ is contained in $\mm$, and the number of geodesics between $M_1$ and $M_2$ in $\mm$ is $m^{m-2}$.

For any other pair $M_1$ and $M_2$ of antipodes in $\mm$, one of the matchings, say $M_1$, has an edge $e$ contained (except for its endpoints) in the interior of $H$. 
Denote by $H_1$ and $H_2$ the two components of $H \setminus e$.
$M_1$ must have at least one edge $e_1$ contained in $H_1$ and another $e_2$ contained in $H_2$. 
Let $e'$ be an edge incident with a vertex of $e_1$ and a vertex of $e_2$ such that the cycle $M_1\cup M_2$ splits in $(M_1*e')\cup M_2$ into two cycles. Then, $M_1*e'$ is a neighbor of $M_1$ in a geodesic between $M_1$ and $M_2$ in $\Pm$. Since $e$ and $e'$ intersect, $M_1*e'\notin \mm$. Thus, the number of geodesics between $M_1$ and $M_2$ in $\mm$ is smaller than $m^{m-2}$.
\end{proof}

\section{Factorization of Permutations by Transpositions}\label{section_cycleFactorization}
Let $S_m$ be the symmetric group on $m$ elements. 
A {\bf minimal factorization by transpositions} (or, simply, \textbf{minimal factorization}) of a permutation $\pi\in S_m$, is a product $\sigma_1\cdots\sigma_n$ of a minimal number of transpositions, such that $\pi=\sigma_1\cdots\sigma_n$.
In this section, the number of geodesics between antipodes in $\Pm$, is shown to be equal to the number of minimal factorizations of an $m$-cycle in $S_m$.
The latter was proved by D\'{e}nes \cite{denes} to be equal to the number of labeled trees on $m$ vertices. 
This result can be deduced from a classic result of Hurwitz \cite{hurwitz}, as explained by Strehl \cite{strehl}.


Throughout this section, the vertices of the underlying complete graph $K_{2m}$ of $\Pm$, are labeled by $\{1,-1,\dots,m,-m\}$.
Denote by $K_{m,m}$ the complete bipartite graph with sides $\{1,\dots,m\}$ and $\{-1,\dots,-m\}$. 
Denote by $\Pmm$ the subgraph of $\Pm$, induced by the perfect matchings of $K_{m,m}$.
Matchings in $\Pmm$ correspond naturally to permutations in $S_m$, by identifying a matching $\pi\in\Pmm$ with a permutation $\pi'\in S_m$ such that $\pi'(i)=j$ for every $(i,-j)\in\pi$.
Using this correspondence, a matching $\pi\in\Pmm$ will also be referred to as a permutation in $S_m$ (and vice versa).

\begin{obs}\label{obs_ins_eq_transpose}
Let $\pi\in\Pmm$. For every $1\leq i\leq j\leq m$,
$$\pi\circ (i,j)=\pi*(i,-\pi(j))=\pi*(j,-\pi(i)),$$
where $(i,i)$ is understood as the identity of $S_m$.
\end{obs}

In Observation \ref{obs_ins_eq_transpose}, the neighbor $\pi*(i,-\pi(j))=\pi*(j,-\pi(i))$ of $\pi$ in $\Pmm$, is one of the two neighbors of $\pi$ in $\Pm$, corresponding to the choice of the two edges $(i,-\pi(i))$ and $(j,-\pi(j))$, as in Remark \ref{rem_regular}. 
The other neighbor in $\Pm$, corresponding to this choice of two edges, is not in $\Pmm$. 
Thus, every pair of edges in $\pi$ corresponds to a unique neighbor of $\pi$ in $\Pmm$, which also corresponds to right multiplication of $\pi$ by a unique transposition in $S_m$.
In other words, adjacency in $\Pmm$ can be understood as right multiplication by transpositions.

Let $G$ be a group, and let $S$ be a symmetric ($S^{-1}=S$) generating set of $G$.
Recall that the (right) {\bf Cayley graph} $X(G,S)$, is the directed graph on the elements of $G$, in which $(g_1,g_2)$ is an edge if $g_1 s= g_2$ for some $s$ in $S$.

\begin{cor}\label{cor_cayley}
$\Pmm$ is the underlying simple graph of the (right) Cayley graph $X(S_m,S)$, where $S$ is the set of all of the transpositions in $S_m$.
\end{cor}

\begin{obs}\label{obs_geodesic_neighbor}
Let $\pi$ and $\sigma$ be two matchings in $\Pmm$, and let $\pi*e$ be a neighbor of $\pi$ (in $\Pm$) in a geodesic between $\pi$ and $\sigma$. Then $\pi*e\in\Pmm$.
\end{obs}
Indeed, the union of $\pi$ and $\sigma$, in Observation \ref{obs_geodesic_neighbor}, is a vertex-disjoint union of alternating cycles of even lengths. 
Clearly, the signs of the labels of the vertices within each cycle alternate. 
By Observation \ref{obs_inc_comp} and Theorem \ref{thm_distance}, a neighbor $\pi*e$ of $\pi$ in $\Pm$, is in some geodesic between $\pi$ and $\sigma$ if and only if a cycle in $\pi\cup \sigma$ is split into two cycles in $\pi*e\cup \sigma$. 
This implies that the vertices of $e$ belong to the same cycle in $\pi\cup \sigma$, and that their labels have different signs (see Figure \ref{fig_split_pmm}).

\begin{figure}[hbt]
\[
\begin{aligned}
\begin{tikzpicture}[scale=0.5]
\fill (2.4,3.4) circle (0.1) node[above]{\tiny -};
\fill (3.4,2.4) circle (0.1) node[right]{\tiny +};
\fill (3.4,1)   circle (0.1) node[right]{\tiny -};
\fill (2.4,0)   circle (0.1) node[below]{\tiny +};
\fill (1,0)     circle (0.1) node[below]{\tiny -};
\fill (0,1)     circle (0.1) node[left]{\tiny +};
\fill (0,2.4)   circle (0.1) node[left]{\tiny -};
\fill (1,3.4)   circle (0.1) node[above]{\tiny +};

\draw (0,1)--(0,2.4);
\draw[blue,dotted] (0,2.4)--(1,3.4);
\draw (1,3.4)--(2.4,3.4);
\draw[blue,dotted] (2.4,3.4)--(3.4,2.4);
\draw (3.4,2.4)--(3.4,1);
\draw[blue,dotted] (3.4,1)--(2.4,0);
\draw (2.4,0)--(1,0);
\draw[blue,dotted] (1,0)--(0,1);

\draw[very thin,densely dashed] (1,3.4)--(1,0);
\node at (1,1.7) {$e$};
\node at (1.7,-1.2) {\small $\pi\cup\sigma$};
\end{tikzpicture}
\end{aligned}
\ \ \ \longrightarrow \ \ \
\begin{aligned}
\begin{tikzpicture}[scale=0.5]
\fill (2.4,3.4) circle (0.1) node[above]{\tiny -};
\fill (3.4,2.4) circle (0.1) node[right]{\tiny +};
\fill (3.4,1)   circle (0.1) node[right]{\tiny -};
\fill (2.4,0)   circle (0.1) node[below]{\tiny +};
\fill (1,0)     circle (0.1) node[below]{\tiny -};
\fill (0,1)     circle (0.1) node[left]{\tiny +};
\fill (0,2.4)   circle (0.1) node[left]{\tiny -};
\fill (1,3.4)   circle (0.1) node[above]{\tiny +};

\draw (0,1)--(0,2.4);
\draw[blue,dotted] (0,2.4)--(1,3.4);
\draw[blue,dotted] (2.4,3.4)--(3.4,2.4);
\draw (3.4,2.4)--(3.4,1);
\draw[blue,dotted] (3.4,1)--(2.4,0);
\draw[blue,dotted] (1,0)--(0,1);

\draw (1,3.4)--(1,0);
\draw (2.4,3.4)--(2.4,0);

\node at (1.7,-1.2) {\small $\pi*e\cup\sigma$};

\end{tikzpicture}
\end{aligned}
\]
\caption{Neighbors in a geodesic in $\Pmm$.}
\label{fig_split_pmm}

\medskip 
\begin{minipage}{0.65\textwidth} 
{\footnotesize Here, the dotted edges belong to $\sigma$ and the solid edges belong to $\pi$.\par}
\end{minipage}

\end{figure}

\begin{cor}\label{cor_cayley_geodesic}
Let $\pi$ and $\sigma$ be two matchings in $\Pmm$. 
Every geodesic between $\pi$ and $\sigma$ in $\Pm$ is contained in $\Pmm$.
\end{cor}

By Corollaries \ref{cor_cayley} and \ref{cor_cayley_geodesic}, the number of geodesics in $\Pm$ between two elements $\pi$ and $\sigma$ in $\Pmm$, is equal to the number of minimal factorizations of $\pi\sigma^{-1}$. In particular, by Theorem \ref{thm_distance}, an $m$-cycle $\pi$ and the identity of $S_m$ are antipodes in $\Pm$, and the number of geodesics between them is equal to the number of minimal factorizations of $\pi$.

\begin{cor}\label{cor_geodesics_factorizations}
The number of geodesics between antipodes in $\Pm$, is equal to the number of minimal factorizations of an $m$-cycle in $S_m$.
\end{cor}

D\'{e}nes' result  was generalized to enumeration of maximal chains in the noncrossing partition lattice of any finite Coxeter type (e.g. \cite{chapoton} Proposition 9).
Appropriate generalizations of Corollary \ref{cor_geodesics_factorizations} to other Coxeter types is most desired.

\section{Acknowledgments}
I would like to thank my supervisors, Ron M. Adin and Yuval Roichman, for their guidance and help.
I would also like to thank Christian Krattenthaler for his many valuable suggestions.


\end{document}